\documentclass{article}
\usepackage[russian]{babel}
\usepackage{amsmath, amssymb, amsfonts,amsthm,comment,graphicx,amscd}
\usepackage[unicode,pdftex]{hyperref}

\newtheorem{theorem}{Теорема}
\theoremstyle{definition}
\newtheorem{problem}{Задача}
\newtheorem{question}{Вопрос}

\title{Корни многочленов и касательные к окружностям}
\author{Андрей Рябичев%
\footnote{\href{https://mailto:ryabichev@179.ru}{ryabichev@179.ru}},
Константин Щербаков}
\date{}

\begin{document}

\maketitle

\section{Правильный треугольник}\label{s:triang}

Рассмотрим многочлен третьей степени от одной переменной $f(x)$, имеющий три корня.
Отметим эти корни на оси абсцисс и проведём через них вертикальные прямые $k_1,k_2,k_3$.
Легко показать, что существует правильный треугольник $A_1A_2A_3$ с вершинами на соответствующих прямых.

Оказывается \cite{mathstodon}, если отметить вписанную окружность $\omega$ треугольника $A_1A_2A_3$, то вертикальные касательные к $\omega$ пройдут через точки экстремума $f(x)$, а вертикальная прямая, проходящая через центр $\omega$, пройдёт через точку перегиба $f(x)$!

\vspace{.4em}

\hfil \includegraphics[width=6cm]{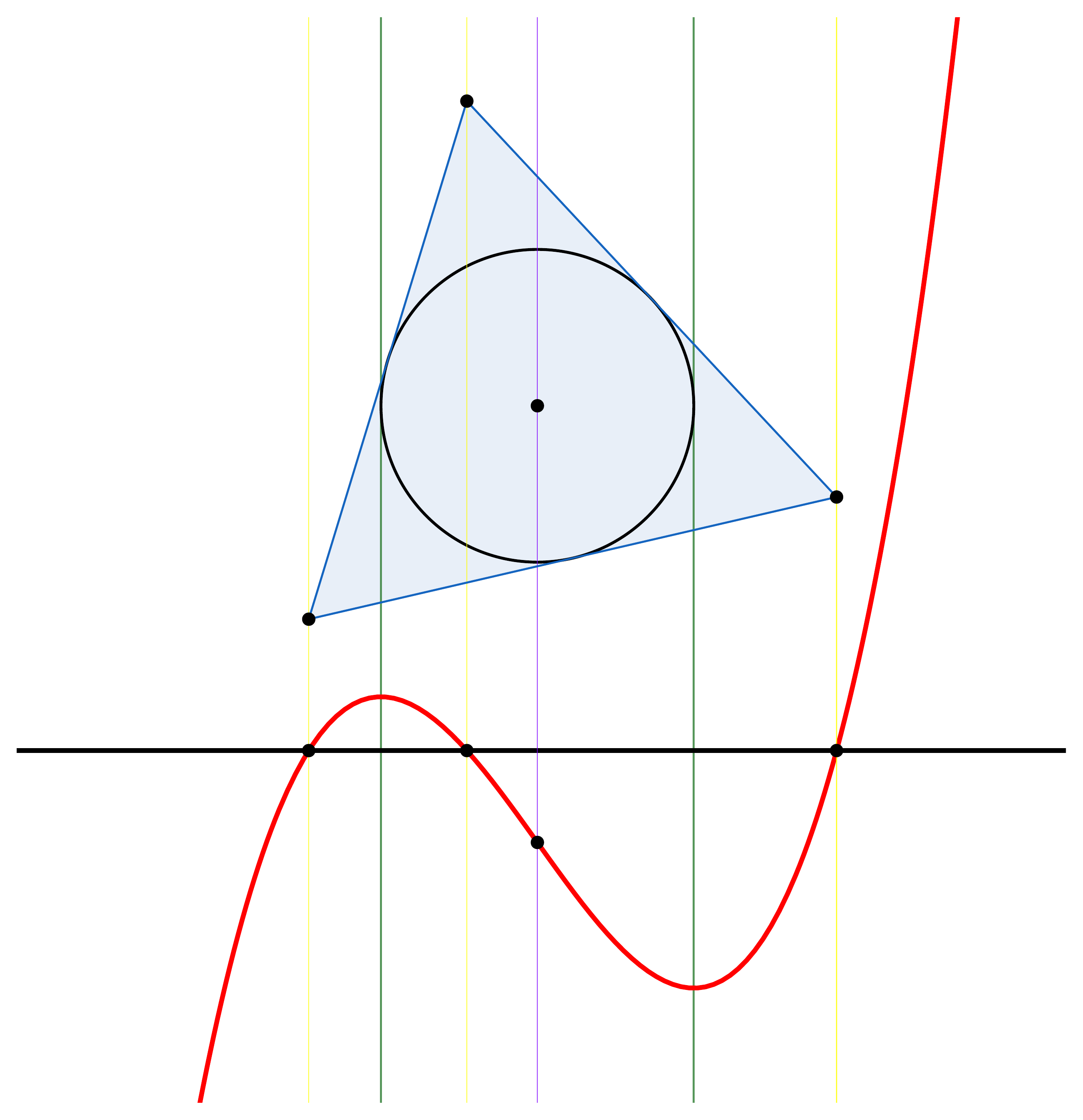}

\vspace{.4em}

Этот удивительный феномен можно пытаться доказать напрямую.
Однако, часто математические утверждения становится легче доказать, если попробовать их обобщить --- нащупать общий ``принцип работы''.
Да и приятнее.
Давайте попробуем сформулировать аналог данного наблюдения для многочленов произвольной степени.

\section{Правильный $n$-угольник}

К сожалению, для многочлена большей степени даже первый шаг построения не проходит:

\begin{problem}
Если $n\ge4$, то не для любых $n$ параллельных прямых существует правильный $n$-угольник с вершинами на этих прямых.
\end{problem}

Проделаем этот шаг наоборот: возьмём правильный $n$-угольник
$A_1\ldots A_n$
и спроецируем его вершины на ось абсцисс.
Пусть $f(x)$ --- многочлен степени $n$, имеющий такой набор корней (некоторые корни могут оказаться кратными; все многочлены с таким условием пропорциональны, нам подойдёт любой из них).

Отметим корни производной $f'(x)$ на оси абсцисс 
и проведём через них вертикальные прямые $l_1,\ldots,l_{n-1}$.
Поскольку корни $f(x)$ имеют кратность не больше $2$, все эти прямые различны.

Следующий факт можно обнаружить, проделав построении в geogebra:

\begin{theorem}\label{3gon}
Самая левая и самая правая прямые из набора $l_1,\ldots,l_{n-1}$ касаются окружности, вписанной в $A_1\ldots A_n$.
\end{theorem}

\vspace{.4em}

\hfil \includegraphics[width=8cm]{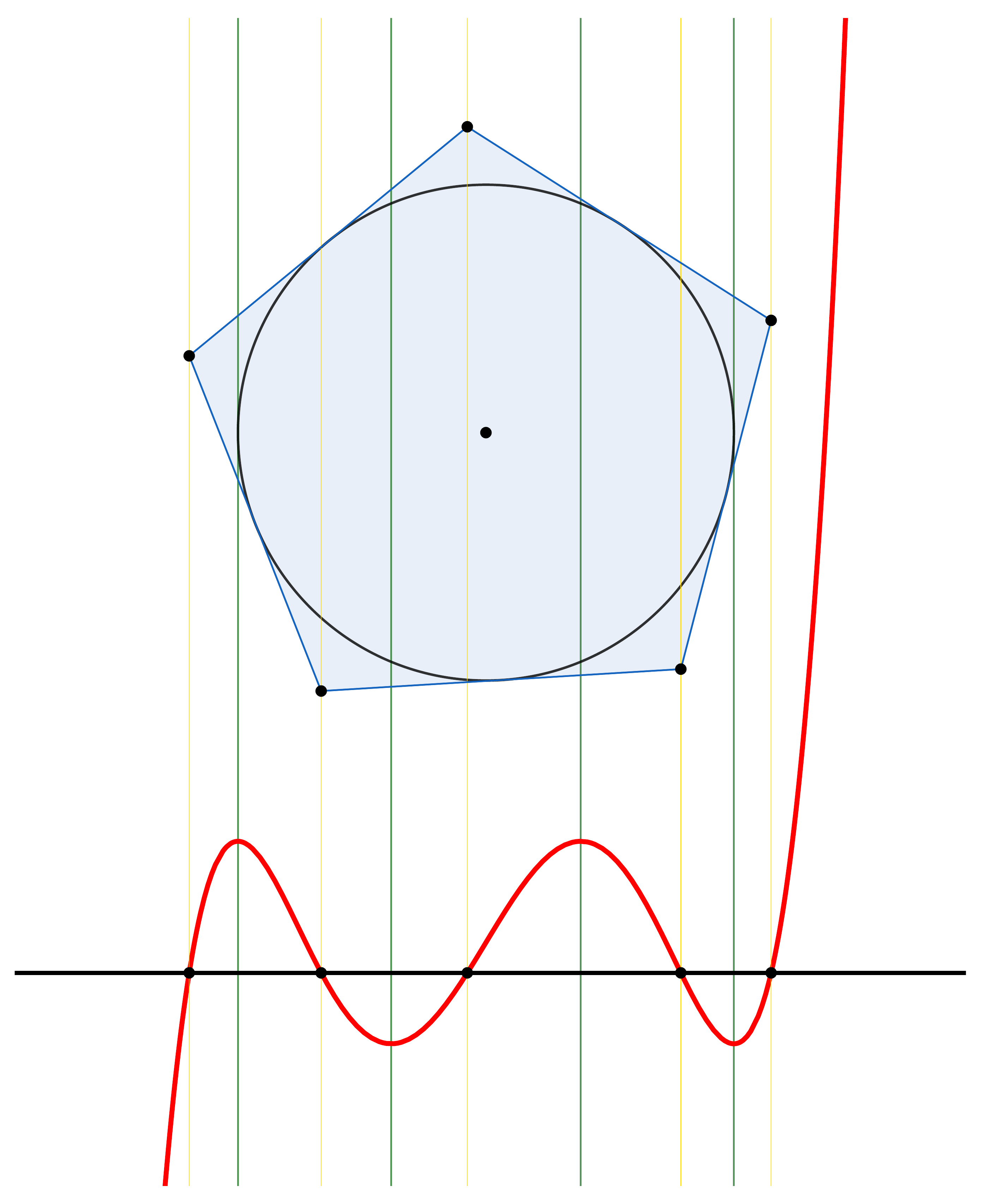}

\vspace{.4em}

Другими словами, если правильный $n$-угольник вращается вокруг своего центра, то корни $f(x)$ будут меняться, но самый большой и самый маленький корни $f'(x)$ остаются неизменными.

\section{Остальные корни производной}

Зададимся вопросом: чему соответствуют остальные корни производной $f'(x)$?
Мысленно вращая $n$-угольник, можно заметить, что $f(x)$ имеет кратные корни ровно в те моменты, когда некоторые диагонали $n$-угольника вертикальны.
Это наблюдение диктует нам следующие построения, приводящие к верной гипотезе.

Отметим все окружности $\omega_1,\omega_2,\ldots$ с центром в центре $O$ многоугольника $A_1\ldots A_n$ и касающиеся какой-нибудь диагонали $A_1\ldots A_n$. Всего мы проведём $\lfloor(n-1)/2\rfloor$ окружностей, включая $\omega$.

\vspace{.4em}

\hfil
\includegraphics[width=8cm]{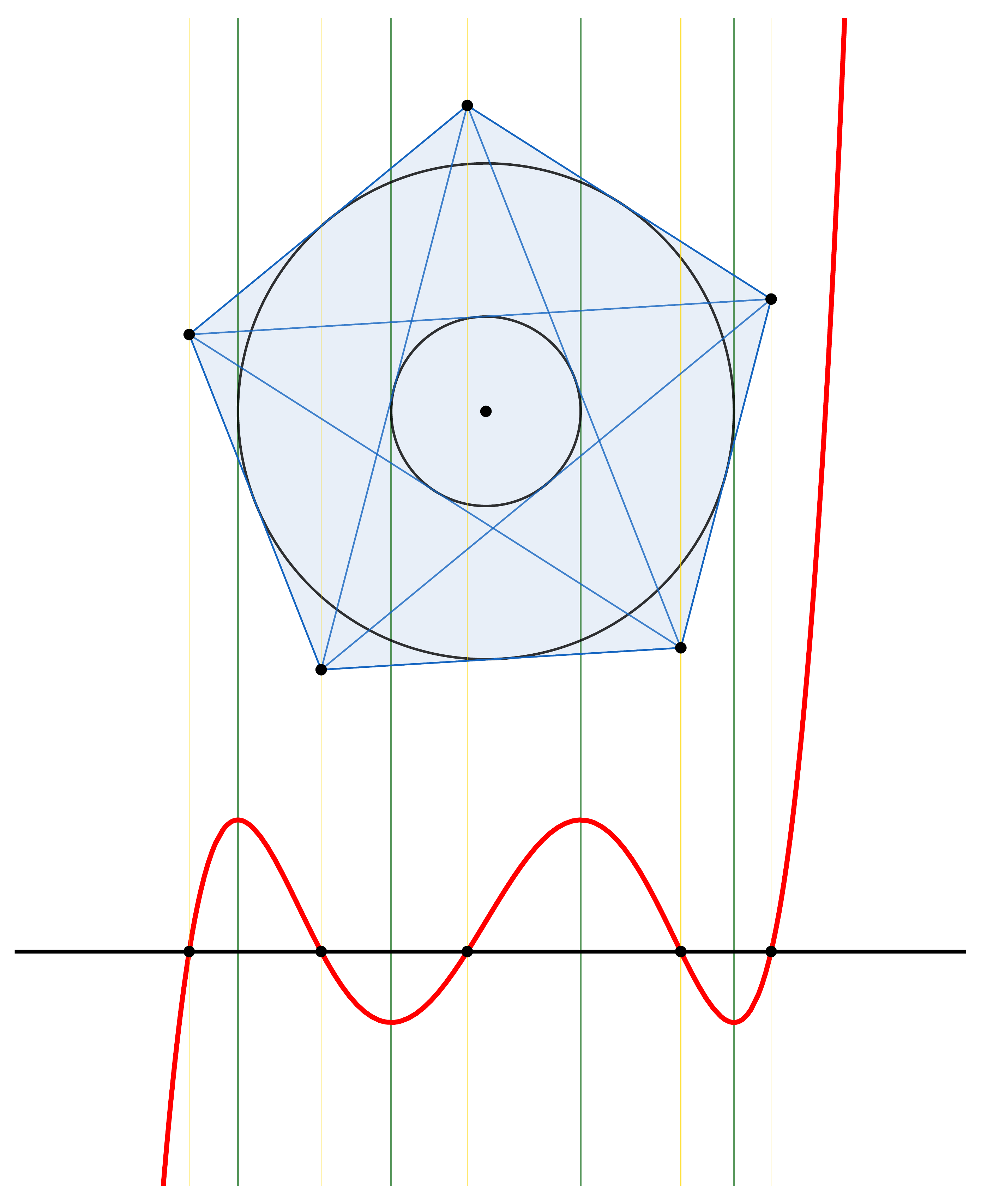}
\vspace{.4em}

\begin{theorem}\label{ngon}
Каждой из окружностей $\omega,\omega_1,\ldots$ касаются две прямые из набора $l_1,\ldots,l_{n-1}$.
\end{theorem}

При чётном $n$ остаётся одна лишняя прямая --- она проходит через $O$ из симметрии.

Отметим, что теорема~\ref{3gon} является частным случаем теоремы~\ref{ngon}.
А она, в свою очередь,
поскольку $f'(x)$ имеет степень $n-1$,
равносильна следующей:

\begin{theorem}\label{deriv}
Будем вращать $A_1\ldots A_n$ вокруг $O$ и в каждый момент строить приведённый многочлен $f(x)$ как выше.
Тогда все коэффициенты $f(x)$ кроме свободного члена не меняются при вращении $A_1\ldots A_n$.
\end{theorem}

Действительно, условие из теоремы~\ref{ngon} однозначно определяет корни $f'(x)$, а значит и сам $f'(x)$, если $f(x)$ приведённый. Из этого следует теорема~\ref{deriv}.

Наоборот, условие из теоремы~\ref{deriv} фиксирует корни $f'(x)$ при вращении многоугольника $A_1\ldots A_n$.
Заметим, что каждая его диагональ в некоторый момент становится вертикальной --- тогда она по построению $f(x)$ проходит через его кратный корень, то есть через экстремум, что влечёт теорему~\ref{ngon}.

\vspace{.5em}

Начнём с прямого доказательства теоремы~\ref{deriv}, требующего лишь умения работать с комплексными числами.

\section{Доказательство через комплексные координаты}\label{s:compl-coord}

Отметим на единичной окружности в $\mathbb{C}$ точки $a_1,\ldots,a_n$ так чтобы они образовывали вершины правильного $n$-угольника,
пронумерованные против часовой стрелки, и $a_1=1$.
Проекция единичного по модулю комплексного числа $a$ на ось абсцисс равна $\mathrm{Re}(a)=\frac12\big(a+\frac1a\big)$.
Таким образом, имеем
$$\textstyle
f(x)=
\left(x-\frac12\big(a_1+\frac1{a_1}\big)\right)
\cdot\ldots\cdot
\left(x-\frac12\big(a_n+\frac1{a_n}\big)\right).
$$

Возьмём комплексное число $t$ по модулю равное единице и умножим все вершины $a_i$ на $t$.
Многочлен примет вид
$$\textstyle
f_t(x)=
\left(x-\frac12\big(ta_1+\frac1{ta_1}\big)\right)
\cdot\ldots\cdot
\left(x-\frac12\big(ta_n+\frac1{ta_n}\big)\right).
$$
Теорема~\ref{deriv} равносильна тому, что все коэффициенты $f_t(x)$ кроме свободного члена не зависят от $t$.
Докажем это.

\begin{problem}\label{xn-2}
Коэффициент при $x^{n-1}$ в $f_t(x)$ равен нулю для всех $t$.
\end{problem}

Рассмотрим коэффициент в $f_t(x)$ при $x^{n-k}$, где $0<k<n$.
Он принимает вид
$$\textstyle
\pm\frac1{2^k}
\sum\limits_{i_1<\ldots<i_k}
\left(ta_{i_1}+\frac1{ta_{i_1}}\right)
\cdot\ldots\cdot
\left(ta_{i_k}+\frac1{ta_{i_k}}\right).
$$
Сгруппируем эту сумму по степеням $t$ (которые могут быть отрицательными, от $-k$ до $k$) и рассмотрим слагаемое, умноженное на $t^m$, $m\ne0$.
Оно будет состоять из всех мономов вида
$$
\dfrac
{a_{i_1}\cdot\ldots\cdot a_{i_j}}
{a_{i_{j+1}}\cdot\ldots\cdot a_{i_k}},
$$
где $m=j-(k-j)$, а индексы выбраны и разбиты на множества $i_1<\ldots<i_j$ и $i_{j+1}<\ldots<i_k$ всеми возможными способами
(число таких способов равно ${n\choose k}\cdot{k\choose j}$).
Обозначим этот коэффициент при $t^m$ через $S$.
Вообще говоря, число $S$ является комплексным.

Пусть $\beta=e^{2\pi i/n}$ --- первообразный корень $n$-й степени из единицы.
Заменим в сумме $S$ каждое $a_i$ на $a_i\beta$.
Очевидно, $S$ при такой подстановке превратится в $S\beta^m$.
Но, с другой стороны, $a_1\beta=a_2$, $a_2\beta=a_3$, \ldots, $a_n\beta=a_1$, поэтому $S$ не изменится.
Следовательно, раз $\beta^m\ne1$, имеем $S=0$.
Это завершает доказательство.
\vspace{.5em}

Использованное нами в конце рассуждение с {\it симметрическими многочленами} довольно стандартно.
Вот другой, более простой пример:

\begin{problem}
Пусть $a_1\ldots,a_n\in\mathbb{C}$ --- вершины правильного $n$-уголь\-ника с центром в нуле.
Докажите, что
$$\textstyle
\text{{\bf (а)} }
\sum\limits_{i_1<\ldots<i_k}
a_{i_1}\cdot\ldots\cdot a_{i_k}=0
\ \text{для любого $k=1,\ldots,n-1$};
$$
$$\textstyle
\text{{\bf (б)} }
\sum\limits_{i\ne j}
a_i^2\cdot a_j=0\
\text{при $n>3$}.
$$
\end{problem}

\section{Многочлены Чебышёва}

Оказывается, у многочлена $f_t(x)$ есть более явный способ его задать, и это сильно упрощает рассуждение.%
\footnote{Авторы благодарят
А.\,Заславского и А.\,Устинова
за это ценное замечание.}

Попробуем задуматься о природе многочлена $f_t(x)$. Обозначим через $\tilde t$ аргумент комплексного числа $t$. По построению $f_t(x)$ обнуляется в $n$ (с учётом кратности) точках вида $x=\cos(\tilde t+2\pi k/n)$ и только в них. Мы строим $f_t(x)$ по его корням, а они все находятся на отрезке $[-1;1]$, поэтому можно положить $x=\cos y$. Тогда $f_t(\cos y)=0$ равносильно равенству $\cos(y)=\cos(\tilde t+2\pi k/n)$, что выполняется ровно тогда же, когда верно $\cos(ny)=\cos(n\tilde t)$. Таким образом, у многочлена $f_t(\cos y)$ корни такие же, как и у выражения $\cos(ny)-\cos(nt)$.


\begin{theorem}\label{th:chebyshov}
    Для каждого $n$ существует многочлен $T_n(x)$ такой, что $T_n(\cos y)=\cos ny$.
\end{theorem}


Эти многочлены называются {\it многочленами Чебышёва первого рода}.
Их существование можно прямо вывести из формулы Муавра:
$(\cos y+i\sin y)^n=\cos(ny)+i\sin(ny)$.
Раскрывая скобки по биному и разделяя выражение на действительную и мнимую часть,
мы видим,
что $\cos(ny)$ выражается как некоторый многочлен от функций $\cos y$ и $\sin y$.
Остаётся заметить, что каждый моном в этом многочлене содержит $\sin y$ в чётной степени, и применить замену $\sin^2y=1-\cos^2y$.

По-другому теорему~\ref{th:chebyshov} можно доказать по индукции, используя следующее соотношение:

\begin{problem}
Для $n\ge1$ верно  $T_{n+1}(x)=2xT_n(x)-T_{n-1}(x)$.
\end{problem}

\begin{problem}
Старший коэффициент в $T_n(x)$ равен $2^{n-1}$.
\end{problem}

Таким образом, многочлены $f_t(x)$ и $\frac1{2^{n-1}}\big(T_n(x)-\cos(n\tilde t)\big)$ степени $n$ оба приведённые и обнуляются в одних и тех же $n$ точках (с учётом кратности). Значит, они совпадают. Из полученного равенства видно, что лишь свободный член $f_t(x)$ зависит от $t$.
Это даёт ещё одно доказательство теоремы~\ref{deriv}.


\vspace{.5em}

Вообще многочлены Чебышёва обладают массой интересных свойств.
Например, $T_n(x)$ имеет наименьшее отклонение от нуля на отрезке $[-1;1]$ среди многочленов степени $n$ с заданным старшим коэффициентом \cite{gashkov}, \cite{tabachnikov}.
В статье \cite{zelevinsky} разобраны много других свойств и соотношений.%
\footnote{
Отдельно отметим связь многочленов Чебышёва с числами Каталана $c_n$, открытую совсем недавно \cite{bychkov}.
Оказывается, если рассмотреть ряд Лорана $\frac{T_m(1/2x)}{T_{m-1}(1/2x)}$,
то при $m\to\infty$ он покоэффициентно приближается к производящей функции $\frac1x+c_0x+c_1x^3+c_2x^5+\ldots$
}
Предлагаем читателю доказать пару из них в качестве упражнения:

\begin{problem}
{\bf (а)}  $D_n(x)=2T_n(x/2)$ является приведённым многочленом с целыми коэффициентами. \ \
{\bf (б)} $D_n(t+\frac{1}{t})=t^n+\frac{1}{t^n}$ для любого $t\neq0$.
\end{problem}

\section{Дальнейшие вопросы и замечания}

Из задачи~\ref{xn-2} можно вывести следующее обобщение утверждения про центр треугольника из \S\ref{s:triang}:

\begin{problem}
Проекция $O$ на ось абсцисс совпадает с корнем $f^{(n-1)}(x)$.
\end{problem}

Также отметим следующий факт, продолжающий это наблюдение:

\begin{problem}
Коэффициенты при $x^{n-3},x^{n-5},\ldots$ в $f(x)$ также равны нулю (кроме свободного члена, если $n$ нечётно).
\end{problem}

Как мы понимаем из теоремы~\ref{deriv}, вторая и последующие производные многочлена $f(x)$ также не меняются при поворотах $n$-угольника $A_1\ldots A_n$.

\begin{question}
Проведём вертикальные прямые $l'_1,\ldots,l'_{n-2}$ через корни $f''(x)$, отмеченные на оси абсцисс.
Как по $A_1\ldots A_n$ построить окружности, касающиеся $l'_1,\ldots,l'_{n-2}$,
не прибегая к функции $f(x)$?
\end{question}

Аналогичные вопросы можно ставить и для остальных производных многочлена $f(x)$ --- и авторам на момент написания статьи неизвестны ответы на них.

\end{document}